\documentclass[a4paper,12pt,reqno]{amsart}
 \pdfoutput=1 
\usepackage{amsmath,amsthm,amssymb}
\usepackage{hyperref,cleveref}
\usepackage{mathrsfs}
\usepackage{graphicx}

\allowdisplaybreaks[1]

\newtheorem{thm}{Theorem}
\newtheorem*{theorem*}{Theorem}

\newtheorem{proposition}[thm]{Proposition}
\newtheorem{remark}[thm]{Remark}

\begin{document}

\title[Wandering domains in $\mathbb{C}\times \mathbb{R}$]{Wandering domains in $\mathbb{C}\times \mathbb{R}$ }
\author[R. Kaur]{Ramanpreet Kaur }
\address{Department of Mathematics, University of Delhi,
Delhi--110 007, India}

\email{preetmaan444@gmail.com}

\begin{abstract}
Approximation theory of entire functions has been used to demonstrate the construction of a map on $\mathbb{C}\times\mathbb{R}$ having wandering domains. We also present suitable modification to this construction that helps in obtaining maps with similar characteristics such as infinitely many attracting Fatou components, infinitely many wandering domains having common paths, etc.
\end{abstract}
\keywords{transcendental entire function, approximation theory, wandering domains,  attracting domains, domains with common paths}
\subjclass[2020]{30D05, 37F10}
\maketitle
\section{Introduction and Preliminaries}
In this note, we shall consider the map
$F:\mathbb{C}\times \mathbb{R}\to\mathbb{C}\times \mathbb{R}$ which is of the following form:
\begin{equation}\label{1}
F(z,w)=(h_1(z,w),h_2(w)),
\end{equation}
where $h_1(z,w)$ is a transcendental entire function in the variables $z,w$ and $h_2(w)$ is a real valued function. For the construction of the required map, first we shall recall few definitions and results. The Fatou set of $F$, denoted by $\mathscr{F}(F)$, is  the largest open set where the family of iterates $\{F^n\}_{n\in\mathbb{N}}$ of $F$ is normal. A connected component of the Fatou set is called a Fatou component. The complement of the Fatou set is called as Julia set. For more details about the Fatou-Julia theory in several variables, one can refer to \cite{Fornaess, Fornaess 2} and the references therein.

In one variable, for a rational function of degree $d\geq 2$, every Fatou component is pre-periodic \cite{Sullivan}. In particular, the periodic Fatou component, say $U$, is one of the following forms.
\begin{enumerate}
\item Attracting: every point of the component converges to periodic point $z_0\in U$ under iteration of the function. In this case, the multiplier at the periodic point has modulus less than 1;
\item Parabolic: every point of the component converges to periodic point $z_0\in \partial U$ under iteration of the function;
\item Siegel disk: the function on $U$ is conformally conjugate to an irrational rotation of a unit disk. In this case, the multiplier at the periodic point  has modulus 1;
\item Herman ring: the function on $U$ is conformally conjugate to an irrational rotation of a standard annulus.
\end{enumerate} 
In contrast to this, a transcendental entire function does not have Herman rings. Indeed they may have 
\begin{enumerate}
\item Baker domains: a Fatou component $U$ for which the orbit of every point of $U$ converges to $z_0$ and the function is not analytic at $z_0$. It is clear that for transcendental entire functions, $z_0=\infty$;
\item Wandering domains: a Fatou component $U$ which is not pre-periodic \cite{Baker1}. 
\end{enumerate}
This led to the study of the classification of Fatou components in higher variables. Several authors have discussed about these components for maps of certain types, for example, recurrent Fatou components of a transcendental H\'enon map \cite{Arosio}, Fatou components of attracting skew products \cite{Peters}, etc. 

In this note, we shall restrict ourselves to wandering domains corresponding to a map $F$ mentioned in \Cref{1}.
 We make use of approximation theory to construct an example of wandering domains in $\mathbb{C}\times\mathbb{R}$. In fact, the construction of  such a map is inspired from the construction of a Carleman set used by Singh \cite{Singh}. We construct a general Carleman set in one variable. This leads us to construct an example of wandering domains in $\mathbb{C}\times \mathbb{R}$. We shall prove the following

\begin{thm}\label{Thm1}
There exists a mapping $F:\mathbb{C}\times \mathbb{R}\to \mathbb{C}\times \mathbb{R}$ having escaping wandering domains.
\end{thm}

To apply approximation theory, let us first recall the definition of a Carleman set (in $\mathbb{C}$). For a relatively closed subset  $G$
of $\mathbb{C}$, and 
\[\mathbb{C}(G)=\{f: G\to\mathbb{C}\mid f \text{ is continuous on } G \text{ and analytic on } G^{\circ}\}.\]
We call $G$ as a Carleman set in $\mathbb{C}$ if for any $f\in \mathbb{C}(G)$ and for any positive, continuous function $\epsilon$ on $G$, there exists an entire function $g$ such that 
\[|f(z)-g(z)|<\epsilon(z) \text { for every } z\in G.\]

The following is a  sufficient condition for a set to be  a Carleman set.
\begin{thm}[Nersesjen, \cite{Gaier}]\label{thm 2}
Suppose $G$ is  a closed proper subset of $\mathbb{C}$. Then $G$ is a Carleman set in $\mathbb{C}$ if and only if $G$ satisfies the following conditions:
\begin{itemize}
\item[(i)] $\hat{\mathbb{C}}\setminus G$ is connected.
\item[(ii)] $\hat{\mathbb{C}}\setminus G$ is locally connected at $\infty$.
\item[(iii)] For every compact set $K\subset \mathbb{C}$, there exists a neighborhood $V$ of $\infty$ in $\hat{\mathbb{C}}$ such that no components of $G^{\circ}$ intersects both $K$ and $V$.
\end{itemize}
\end{thm}
\section{Preparation for the Main Result}
To construct the desired map, first we construct an infinite number of pairwise disjoint domains with certain properties, so that no two domains are subsets of a single Fatou component of $F$. We prove that our desired map $F$ with wandering domains have the following form:
\[F(z,w)=(e^{f_1(z)}+\delta e^{\text{Re}f_2(w)}, e^{\text{Re}f_2(w)}), \text{ for every }(z,w)\in\mathbb{C}\times \mathbb{R}.\] 
Here, $f_1, f_2$ are entire functions (existence of such functions is shown during the construction of a Carleman set) and $\delta\in(0,1)$ is chosen suitably (one of the conditions is $\frac{17\delta}{16}<\frac{9}{8}$).\\
\subsection*{Construction of a Carleman set in one variable:} 
 For every $k\in\mathbb{N}$, define 
\[ B_k:=\{z\in\mathbb{C}: |z-\delta_k|\leq l_k\},\]
where the sequences $\{\delta_k\}$ and $\{l_k\}$ are chosen such that
\begin{itemize}
\item $\{\delta_k\}$ and $\{l_k\}$ are strictly  increasing real sequences converging to infinity.
\item  $l_1$ satisfies $\frac{17\delta}{16}<\frac{l_1}{8}$.
\item $B_i\cap B_j=\emptyset$ for every $i\not=j$.
\item $B_1\cap G_1=\emptyset$,  where  $G_1=\{z\in\mathbb{C}: |z-1|\leq 3\}$.
\end{itemize} 
Now, choose a sequence $\{m_k\}$ such that for every $k\in\mathbb{N}$, we have
\[\delta_k+l_k<m_k<\delta_{k+1}-l_{k+1}.\]
For such a sequence, define 
\[M_k:=\{z\in\mathbb{C}: \text{Re}z= m_k\}, k\in\mathbb{N}.\]
Further, define
\begin{align*}
A_k&:=\{z\in\mathbb{C}: |z+r_k|\leq s_k\}, k\in\mathbb{N}\\
 \text{ and, }L_k&:=\{z\in\mathbb{C}: \text{Re} z= -t_k\}, k\in\mathbb{N}
\end{align*}
where, 
\begin{itemize}
\item $\{-r_k\}$ is chosen to be a decreasing real sequence converging to $-\infty$.
\item $\{s_k\}$ is chosen such that $A_i\cap A_j=\emptyset$ for every $i\not=j$.
\item $\{t_k\}$ is chosen such that $(-r_{k+1}+s_{k+1})<-t_k<(-r_k+s_k)$, for every $k\in\mathbb{N}$.
\end{itemize}
Using the continuity of the exponential function $z\mapsto e^z$, we have 
\begin{itemize}
\item[(a)] for  every $\epsilon_k= \frac{l_{k+1}}{4}$, there exists $\delta'_k>0$ such that 
\[|z-\log (\delta_{k+1})|< \delta'_k \text{ implies that } |e^z- \delta_{k+1}|<\frac{l_{k+1}}{4}.\]
\item[(b)] As $\frac{s_1}{8}>0$, there exists $\delta^{''}_1>0$ such that 
\[|z-\pi\iota+\log (r_{1})|< \delta^{''}_1 \text{ implies that } |e^z+ r_{1}|<\frac{s_{1}}{8}.\]
\item[(c)] there exists $\delta^{''}_2>0$ such that 
\[|z|<\delta^{''}_2 \text{ implies that }|e^z-1|<\frac{1}{16}.\]

\end{itemize}
\textbf{Defining the functions $g(z)$ and $\epsilon(z)$:}
Let,
\begin{equation*}
g(z)=\left\{
\begin{array}{ll}
\log (\delta_{k+1}), & z\in B_k,k\in\mathbb{N}\\
 \pi \iota+\log (r_1),& z\in A_k, k\in\mathbb{N}\\
 \log 1, &z\in G_1\cup\left(\cup_{k\in\mathbb{N}} (M_k\cup L_k) \right),
\end{array}
\right.
\end{equation*}
and
\begin{equation*}
\epsilon(z)= \left\{
\begin{array}{ll}
\delta'_k, &z\in B_k,k\in\mathbb{N}\\
\delta^{''}_1,& z\in A_k, k\in\mathbb{N}\\
\delta^{''}_2,& z\in G_1\cup\left( \cup_{k\in\mathbb{N}}(M_k\cup L_k)\right). \end{array}
\right.
\end{equation*}  
Further, let $G= G_1\cup \left(\cup_{k\in\mathbb{N}}(B_k\cup M_k\cup A_k\cup L_k)\right)$ (See Figure 1). By \Cref{thm 2}, $G$ is a Carleman set. Since $g\in\mathbb{C}(G)$, there exists an entire function $f_1(z)$ such that 
\[|f_1(z)-g(z)|<\epsilon(z),\text{ for every }  z\in G.\]
In particular, we get that 
\begin{itemize}
\item[(d)] $|f_1(z)-\log (\delta_{k+1})|< \delta'_k$, whenever $z\in B_k, k\in\mathbb{N}$.
\item[(e)] $|f_1(z)-\pi \iota+\log (r_1)|<\delta^{''}_1$, whenever $z\in A_k, k\in\mathbb{N}.$
\item[(f)] $|f_1(z)-1|<\delta^{''}_2$, whenever $z\in G_1\cup\left( \cup_{k\in\mathbb{N}}(M_k\cup L_k)\right)$.
\end{itemize}   
 This further implies that  
 \begin{itemize}
\item[(g)] $|e^{f_1(z)}- \delta_{k+1}|<\frac{l_{k+1}}{4}$, whenever $z\in B_k, k\in\mathbb{N}$.
\item[(h)] $|e^{f_1(z)}+r_1|< \frac{s_1}{8}$, for every $z\in A_k, k\in\mathbb{N}$.
\item[(i)] $|e^{f_1(z)}-1|<\frac{1}{16},$ for $z\in G_1\cup (\cup_{k\in\mathbb{N}}( M_k\cup L_k))$.
\end{itemize}
Now, using continuity of $e^w$ at origin, there exists $d_1>0$ such that whenever $|w-\log 1|<d_1$, we have $|e^w-1|<\frac{1}{16}.$\\
Further, define 
$g_1(w)=\log 1$ and $\epsilon_1(w)=d_1,\text{ for } w\in\mathbb{R}$. Since $\mathbb{R}$ is a Carleman set \cite{Gaier}, there exists an entire function $f_2$ such that 
\[|f_2(w)-g_1(w)|<\epsilon_1(w), \text{ for every }  w\in\mathbb{R}.\] 
In particular, 
\[|f_2(w)-\log 1|<d_1, \text{ for every }   w\in\mathbb{R}.\]
This further implies that 
\[|e^{f_2(w)}-1|<\frac{1}{16}, \text{ for every }   w\in\mathbb{R}.\]
From here, we can easily observe that  
\[|e^{\text{Re}f_2(w)}|=|e^{f_2(w)}|<\frac{17}{16}, \text{ for every } w\in\mathbb{R}.\]
\begin{figure}
\centering
\includegraphics[totalheight=8cm]{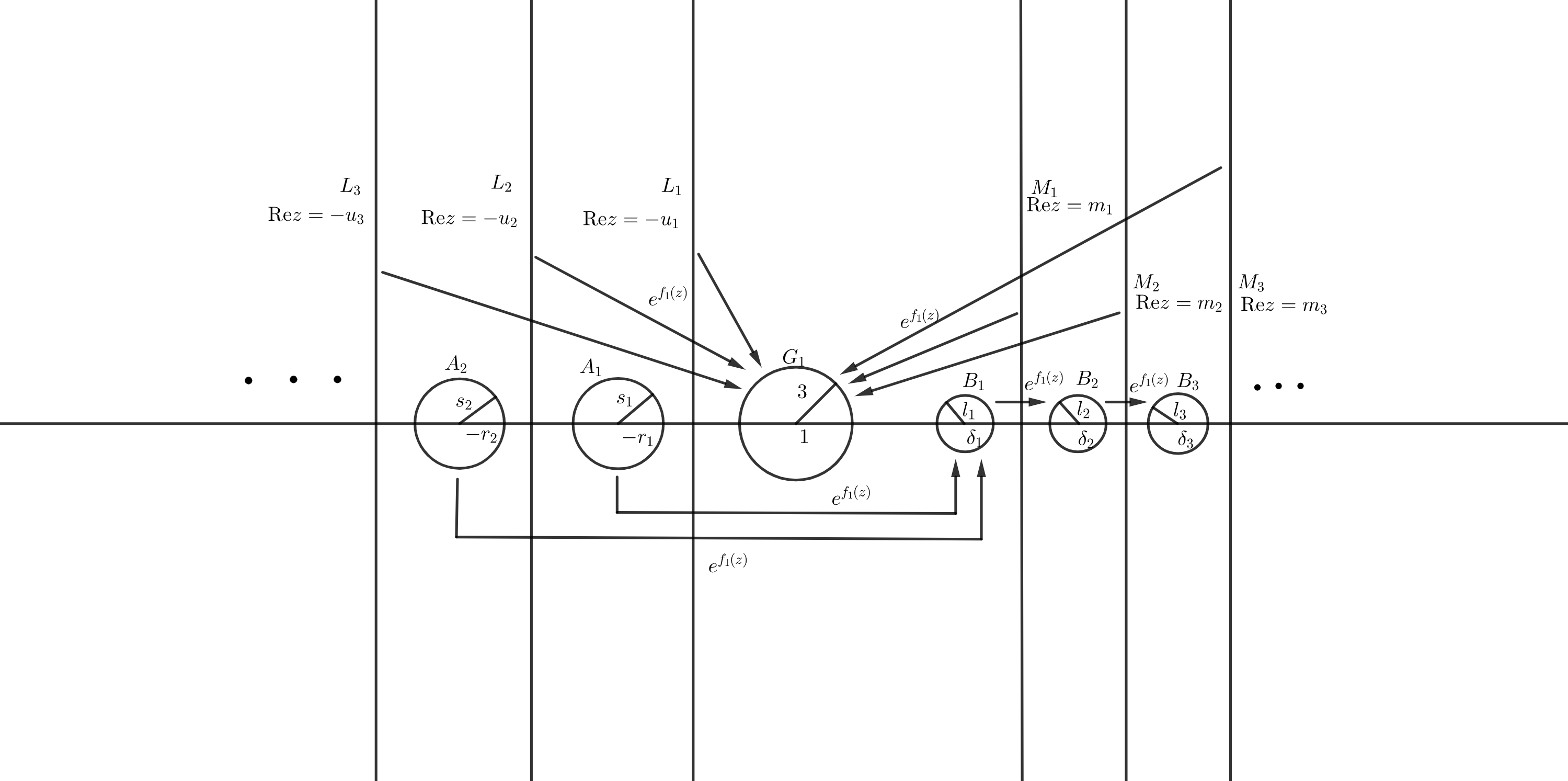}
\caption{ The arrows represent mapping of $L_k, M_k, B_k$ (in the proof of \Cref{Thm1}) and that of $A_k$ (in the proof of \Cref{prop}) under the function $e^{f_1(z)}$, for $k\in \mathbb{N}$.
}
\label{fig:1}
\end{figure}

\begin{proof}[Proof of \Cref{Thm1}]
 
Let $f_1$ and $f_2$ be defined as above. Let 
\[F(z,w)= (e^{f_1(z)}+\delta e^{\text{Re}f_2(w)}, e^{\text{Re}f_2(w)}), \text{ for } (z,w)\in\mathbb{C}\times\mathbb{R}.\]
For wandering domains of the above  map, define  
\begin{align*}
B_k'&:=\{(z,w)\in\mathbb{C}\times\mathbb{R}: |z-\delta_k|\leq l_k, |w|<\frac{l_k}{8\delta }\}, k\in\mathbb{N}\\
M_k'&:=\{(z,w)\in\mathbb{C}\times\mathbb{R}: \text{Re}z=m_k , w\in\mathbb{R}\}, k\in\mathbb{N}\\
A_k'&:=\{(z,w)\in\mathbb{C}\times\mathbb{R}: |z+ r_k|\leq s_k, |w|\leq \frac{l_1}{4}\}, k\in\mathbb{N}\\
L_k'&:=\{(z,w)\in\mathbb{C}\times\mathbb{R}: \text{Re}z=-t_k, w\in\mathbb{R} \}, k\in\mathbb{N}\\
\text{and } G_2&:=\{(z,w)\in\mathbb{C}\times\mathbb{R}: |z-1|\leq 3, |w|\leq \frac{9}{8\delta}\}.
\end{align*}
First, we shall show that for every $(z,w)\in B_k'$, we have $F(z,w)\in B_{k+1}'$. For this, suppose that $(z,w)\in B_k'$. By definition of $F$, we have $F(z,w)= (e^{f_1(z)}+\delta w, e^{f_2(w)})$. Consider, 
\begin{align*}
|e^{f_1(z)}+\delta e^{\text{Re}f_2(w)}-\delta_{k+1}|&\leq |e^{f_1(z)}-\delta_{k+1}+\delta |e^{\text{Re}f_2(w)}|\\
&\leq \frac{l_{k+1}}{4}+ \delta\frac{17}{16}\\
&\leq \frac{l_{k+1}}{4}+\frac{l_{k+1}}{8}\\
&\leq \frac{l_{k+1}}{2}.
\end{align*}
Also,
\begin{align*}
|e^{\text{Re}f_2(w)}|&\leq \frac{17}{16}\\
&\leq \frac{l_{k+1}}{4\delta}.
\end{align*}
This shows that $F(B_k')\subset B_{k+1}'$, for every $k\in\mathbb{N}$, i.e., $F$ maps $B_k'$ into a smaller poly disc inside $B_{k+1}'$, for every $k\in\mathbb{N}$. Now, assume that $(z,w)\in G_2\cup(\cup_{k\in\mathbb{N}} (M_k'\cup L_k'))$. Hence, we have 
\begin{align*}
|e^{f_1(z)}+\delta e^{\text{Re}f_2(w)}-1|&\leq |e^{f_1(z)}-1|+\delta |e^{\text{Re}f_2(w)}|\\
&\leq \frac{1}{16}+ \delta \frac{17}{16}\\
&<3.  
\end{align*}
Also, 
\begin{align*}
|e^{\text{Re}f_2(w)}|&< \frac{17}{16}\\
&<\frac{9}{8\delta}.
\end{align*}
This, in particular shows that $F(M_k')\subset\{(z,w)\in\mathbb{C}^2: |z-1|< 2, |w|<\frac{9}{8\delta}\}\subset G_2$, for every $k\in\mathbb{N}$.
Since, each $B_k'$ is disjoint and $F(M_k')\subset G_2$, for every $k\in\mathbb{N}$. Let $N_k$ be the Fatou component of $F$ containing $B_k'$, for $k\in\mathbb{N}$. By the construction of $F$, it is clear that each $N_k$ is an escaping wandering domain of $F$ for every $k\in\mathbb{N}$.
\end{proof}

\begin{proposition}
 There exists a mapping $F$ on $\mathbb{C}\times \mathbb{R}$ having infinitely many attracting Fatou components.
\begin{proof} Consider $B_k$ and $M_k$ for every $k\in\mathbb{N}$ as we defined earlier. Now we shall follow the following steps:\\
Step1: Using continuity of $e^z$, 
 there exists $h_k'>0$ such that 
\[|z-\log (\delta_k)|<h_k' \text{ implies that } |
e^{z}-\delta_k|<\frac{l_k}{4}.\]
Step2: Define functions $g(z)$ and $\epsilon(z)$ as follows:
\begin{equation*}
g(z) = 
\left\{
    \begin{array}{ll}
        \log (\delta_k), & z\in B_k, k\in\mathbb{N}\\
        \log 1, & z\in G_1\cup (\cup_{k\in\mathbb{N}}M_k),
        \end{array}
\right. 
\end{equation*}
and
\begin{equation*}
\epsilon(z)=
\left\{
\begin{array}{ll}
\frac{l_k}{4}, &z\in B_k, k\in\mathbb{N}\\
\delta^{''}_2, & z\in G_1\cup (\cup_{k\in\mathbb{N}}M_k).
\end{array}
\right.
\end{equation*}
For such functions, using the definition of a Carleman set,  we have
\[|f_1(z)-\log (\delta_{k})|< \delta'_k, \text{ whenever }  z\in B_k, k\in\mathbb{N}.
\] This further implies that  
\[|e^{f_1(z)}- \delta_{k}|<\frac{l_{k}}{4}.\] 
Similarly, for $z\in G_1\cup(\cup_{k\in\mathbb{N}} M_k)$, we have  
\[|e^{f_1(z)}-1|<3.\]
Now, define $B_k'$ and $M_k'$ for every $k\in\mathbb{N}$ as we did earlier. Consider, 
\begin{align*}
|e^{f_1(z)}+\delta e^{\text{Re}f_2(w)}-\delta_k|&\leq |e^{f_1(z)}-\delta_k|+\delta|e^{\text{Re}f_2(w)}|\\
&\leq \frac{l_k}{4}+\frac{17\delta}{16}\\
&\leq \frac{l_k}{4}+\frac{l_k}{4}\\
&=\frac{l_k}{2},
\end{align*}
and
\begin{align*}
 |e^{\text{Re}f_2(w)}|&\leq \frac{17}{16}\\
 &< \frac{l_k}{8\delta}.
 \end{align*} 
These inequalities show that $F(B_k')\subset B_k'$ for every $k\in\mathbb{N}$, i.e., every $B_k'$ is being mapped by $F$ into a smaller poly disc inside $B_k'$. For $k\in\mathbb{N}$, let $C_k$  be the Fatou component containing $B_k'$ which is an attracting Fatou component. Since, $F(G_2\cup(\cup_{k\in\mathbb{N}}M_k'))\subset G_2$ and each $B_k'$ is disjoint, we have $C_i\cap C_j=\emptyset$ for $i\not=j$. This gives us the existence of an infinite number of  attracting Fatou components.
\end{proof}
\end{proposition}
Recall that two wandering domains $C_0$ and $V_0$ are said to have a common path if $C_m=V_n$ for some $n,m\in \mathbb{N}$, where $C_m $ and $V_n$ are the Fatou components containing $F^m(C_0)$ and $F^n(V_0)$ respectively. If $C_m\not=V_n$ for any $n,m\in \mathbb{N}$, then $C_0$ and $V_0$ are said to have distinct paths.
\begin{proposition}\label{prop}
There exists a mapping  $F$ on $\mathbb{C}\times \mathbb{R}$ with infinitely many wandering domains having a common path.
\begin{proof}
For this, we shall make use of the domains $A_k', k\in\mathbb{N}$.
Firstly, using continuity  of $e^z$ there exists $\delta^{''}_3>0$ such that 
\[|z-\log (\delta_1)|< \delta^{''}_3 \text{ implies that } |e^z-\delta_1|<\frac{l_1}{4}.\]
Define the function $g(z)$ and $\epsilon(z)$ as 
\begin{equation*}
g(z)=
\left\{
\begin{array}{ll}
\log(\delta_1), & z\in A_k,k\in\mathbb{N}\\
\log (\delta_{k+1}), & z\in B_k,k\in\mathbb{N}\\
\log 1, &z\in G_1\cup\left(\cup_{k\in\mathbb{N}} (M_k\cup L_k) \right),
\end{array}
\right.
\end{equation*} 
and 
\begin{equation*}
\epsilon(z)= 
\left\{
\begin{array}{ll}
\frac{l_1}{4}, &z\in A_k, k\in\mathbb{N}\\
\delta'_k, &z\in B_k,k\in\mathbb{N}\\
\frac{1}{4},& z\in G_1\cup\left( \cup_{k\in\mathbb{N}}(M_k\cup L_k)\right).
\end{array}
\right.
\end{equation*}  
Now, for every $k\in\mathbb{N}$, define
\[A_k':=\{(z,w)\in\mathbb{C}\times \mathbb{R}: |z-r_k|\leq s_k, |w|\leq \frac{l_1}{4}\}.\]
Arguing as in the previous manner, we can observe that
$|e^{f_1(z)}+\delta e^{f_2(z)}-\delta_1|<\frac{l_1}{2}$ and $|e^{f_2(w)}|<\frac{l_1}{4}$, for every $(z,w)\in A_k'$. This means that $F(A_k')\subset B_1'$, for every $k\in\mathbb{N}$. In particular, if $D_k$ is the Fatou component containing $A_k'$,  $k\in\mathbb{N}$, then $C_1$ (the Fatou component containing $B_1'$) is a common path for these  wandering domains. 
\end{proof}
\end{proposition}
\begin{remark} From the above proposition, one can observe that by modifying functions $g(z), \epsilon(z)$ and by suitably using domains, we can have
\begin{itemize}
\item An $F$ with wandering domains each having a distinct path.
\item $F$ having wandering domains with a common path as well as having wandering domains with distinct paths.
\end{itemize}
\end{remark}
\section*{Acknowledgement}
I would like to thank my thesis advisor, Prof. Sanjay Kumar, for several fruitful discussions. 
This work is supported by the National Board for Higher Mathematics, India.

\end{document}